\documentclass[11pt]{amsart}
\usepackage{amsmath,amsfonts,amssymb,amscd,amsthm,amsbsy,epsf}
\newtheorem{theorem}{Theorem}[section]

\theoremstyle{remark}                  

\def\IH{{\mathbb H}}
\def\IR{{\mathbb R}}
\def\R{\IR}
\def\YMH{{\mathcal Y \mathcal M \mathcal H}}
\def\SYMH{{\mathcal S\mathcal Y\mathcal M\mathcal H}}
\def\overex{\vec x}

\def\overdex {\vec {dx}}
\def\oversigma{\vec \sigma}
\def\sigmadotex{\vec \sigma \cdot \vec x}
\def\sigmadotdex{\vec \sigma \cdot \vec {dx}} 
\def\excrosssigma{\vec x \times \vec \sigma}
\def\excrossdex {\vec x \times \vec{dx}}

\begin{document}
\title{Spherically symmetric solutions of a boundary value 
	problem for monopoles}
\author{Antonella Marini}
\address{Dipartimento Di Matematica,	
Universita' di L'Aquila,
67100 L'Aquila, Italy}
\email{marini@univaq.it}
\author{Lorenzo Sadun}
\address{Department of Mathematics,
University of Texas,
Austin, TX}
\email{sadun@math.utexas.edu}
\begin{abstract}
In this paper we study spherically symmetric monopoles, which are critical
points for the Yang-Mills-Higgs functional over a disk
in 3 dimensions, with prescribed degree and covariant constant at the boundary.
This is a 3-dimensional gauge-theory generalization of the Ginzburg-Landau 
model in 2 dimensions. 
\end{abstract}

\maketitle
\pagestyle{headings}            
\markboth{A. Marini, L. Sadun}{Yang-Mills-Higgs, symmetric bundles} 

\baselineskip=18pt		
\section{Introduction}
\setcounter{equation}{0}
In this paper we treat a 3-dimensional analogue of the vortex equations in 
2 dimensions and look for solutions with spherical symmetry as described in 
\S2.  The domain considered is a 3-dimensional disk and we prescribe the 
degree of the monopole at the boundary.  Unlike in 2 dimensions, where
Ginzburg-Landau-type functionals appear either with or without gauge
potentials, the problem in 3 dimensions is well-posed only in gauge
theory.  Without a curvature term in the action, a minimizing sequence
for the action would yield a trivial limit \cite{mar1}.  In the
presence of a gauge potential, this problem is well-posed, natural and
has physical meaning.

The most general Yang-Mills-Higgs functional takes the form
\begin{equation}
\label{ymhgeneral}
\YMH(A,\phi) = \frac{\epsilon }{2}\| F\|^2_{L^2} + 
\frac{\rho}{2}\| D_A\phi\|^2_{L^2} +
\frac{\lambda}{8} \|\,|\phi|^2-a^2\|^2_{L^2}\ ,
\end{equation}
for appropriate constants $\epsilon$, $\rho$, $\lambda$ and
$a$. Working on $\R^3$, one usually applies a rescaling of $\phi$, a
rescaling of space, and a rescaling of the action to set
$\epsilon=\rho=a=1$, so the action functional depends on a single
parameter, $\lambda$.  On the unit ball, however, we cannot rescale
space, so we can only eliminate two of the four parameters. We 
set $\rho=a=1$, and obtain a 2-paramater family of functionals
\begin{equation}
\label{ymh}
\YMH_{\epsilon,\lambda} (A,\,\phi)=
\frac{\epsilon }{2}\| F\|^2_{L^2(B^{3})} + \frac{1}{2}\| D_A\phi\|^2_{L^2(B^{3})} +
\frac{\lambda}{8} \|\,|\phi|^2-1\|^2_{L^2(B^{3})}\ .
\end{equation}
(Alternatively, we could work on a sphere of radius $R$.  One can then
rescale to set $\epsilon=1$, at the cost of varying $R$.  We then obtain
a 2-parameter family of functionals indexed by $\lambda$ and $R$.)
We know from the general theory for monopoles (cf. \cite {mar1} for
$\epsilon = 1,$ $\lambda \geq 0$) that there exists a minimum for this
functional which satisfies the Euler Lagrange equations 
\begin{equation} \label{pde}
    \begin{split}
\epsilon \, *D_A* F&=[D_A\phi,\,\phi]\\
*D_A*D_A\phi&=\frac{\lambda}{2}(|\phi|^2-1)\phi
\end{split}
\end{equation} 
and suitable boundary conditions on $\partial B^{3}\equiv S^{2}$ 
(cf. Section 4 and \cite{mar1, mar2}), and is smooth.
In this paper we prove the existence, and describe the form, of 
spherically symmetric solutions to these equations. 

We note that, even for $\lambda=0$, these are not 
solutions to the Bogomolnyi equations found in \cite{ja-ta}.
The Bogomolnyi solutions are obtained only in the limit
$\lambda \to 0$, $R \to \infty$, or equivalently
$\lambda \to 0$, $\epsilon \to 0$.

\section{Spherically symmetric connections, monopoles, and gauge 
transformations}
\setcounter{equation}{0}

We work on the trivial principal $SU(2)$-bundle $P =
B^{3} \times SU(2)$ and its associated vector bundles.  $A$ is a
connection on $P$, which can be viewed as a 1-form on $B^3$ with
values in $su(2)$.
%
%
The Higgs field $\phi$ is a section of the adjoint bundle, i.e., a map
$\phi\,:\, B^{3}\to su(2)$.  Here $su(2)\equiv \{X\in \mathcal
M_{2\times 2}\; :\; tr\, X= 0\,; \, X+{\bar X}^{T}= 0\}$ is the Lie
algebra of $SU(2)$. We identify $su(2)$ with the imaginary quaternions
$Im \,\IH \equiv
\{ x_1i + x_2 j + x_3 k : \,(x_1,x_2,x_3) \in \IR^3\}$ 
by identifying the matrices 
$\begin{pmatrix} 0&1\\ -1&0\end{pmatrix}$, 
$\begin{pmatrix} 0&i\\ i&0\end{pmatrix}$ 
$\begin{pmatrix} i&0\\ 0&-i\end{pmatrix}$,
with $i,j,k$, respectively, and extending this mapping to a
Lie-algebra isomorphism.  It is also convenient to define the 3-vector
of Lie-algebra elements $ \vec \sigma = (i,j,k)$.

The symmetry group $SO(3)$ acts on pairs ($A$, $\phi$), simultaneously
rotating the 3-dimensional base space and the 3-dimensional Lie
algebra.  That is, the triple $(i,j,k)$ transform in the same way as
the triple $(x_1, x_2, x_3)$, so quantities such as $\vec \sigma \cdot
\vec x$ are invariant.

We are interested in finding Yang-Mills-Higgs fields $(A,\phi)$ which
are invariant under this group action.  To this purpose, one needs to
specify the value of the connection one-form \\ $A\,:\, B^{3}\to
\Lambda^{1}(B^{3})\otimes Im \,\IH $ and of the Higgs field $\phi \,:\,
B^{3}\to Im \,\IH$ at one point of each group orbit (on the base) and
impose invariance under the isotropy group of that point.  We find it
convenient to fix the values of $(A,\phi)$ on the slice
$L\equiv\{\vec x \in B^{3}\,:\,
x_{2}=x_{3}=0, x_{1}\in (0,1]\}.$ The isotropy group at $(x, 0, 0)\in
L$ is $SO(2),$ i.e. rotations about the $x_{1}$-axis (and about the $i$ axis
in the Lie-algebra).  For the
Higgs field $\phi$ one has in general
$$\phi (x, 0, 0)= \varphi_{1}(x) i + \varphi_{2}(x) j + \varphi_{3}(x) k\in Im 
\,\IH 
\simeq su(2)\,,$$
where $\varphi_{l}(x)\,,$ $l=1,2,3$ are real-valued functions.
Imposing invariance under $SO(2)$ forces 
$ \varphi_{2}(x) =  \varphi _{3}(x) =0$ for all $x\in (0, 1].$
Applying the action of $SO(3),$ one obtains the symmetric form of 
the Higgs field $\phi$ 
\begin{equation}
    \label{symphi}
\phi = \frac {\varphi(r)}{r} \sigmadotex \equiv 
\frac{\varphi (r)}{r}\, (x_1i + x_2 j + x_3 k)\,,
\end{equation}
with $r\equiv \vert \vec x \vert,$ $\vec x\equiv 
(x_{1}, x_{2}, x_{3})\in \IR^{3},$ $\vec \sigma\equiv (i, j, k)\in Im \IH.$

An $su(2)$-valued connection on the slice 
$L$ is given in general by 
\begin{equation}
    \begin{split}
A(x, 0, 0) &= a_{11}(x) i\, dx_{1} + a_{12}(x) i\, dx_{2}+a_{13}(x) 
i\, dx_{3}
+a_{21}(x) j\, dx_{1}+\\
&+a_{22}(x) k\, dx_{2} 
+a_{23}(x) j \,dx_{3}+a_{31}(x) k \,dx_{1} + a_{32}(x) k\, dx_{2}+
+a_{33}(x) k\, dx_{3} 
\end{split}
\end{equation}
Imposing SO(2)-invariance yields 
	     $a_{12}=a_{13}=a_{21}=a_{31}=0,$ $a_{22}= a_{33},$ and
	     $a_{23}= - a_{32}$.
	     Thus the final form of the $SO(2)-$invariant connection 
	     evaluated at the points $(x,0,0)\in L$ is 
$$A(x, 0, 0)= a(x) i \,dx_1 + b(x) (j\,dx_2 + k\,dx_3)
+c(x) (k\,dx_2  -j\,dx_3) \,.$$
Transporting this slice via the $SO(3)$ group action on $B^{3}$ one 
obtains the invariant version 
\begin{equation}
 A( x_{1}, x_{2}, x_{3})= \frac{\alpha(r)}{r} \vec \sigma \cdot \vec {dx} +
 \frac{\beta(r)}{r^{3}}(\vec x \times \vec {dx})\cdot 
(\vec x \times \vec \sigma)
 +\frac{\gamma(r)}{r^{2}} \vec \sigma \cdot ( \vec x \times \vec {dx})\,,  
 \end{equation}
where``$\times$'' denotes the cross product of vectors and  
$\vec {dx} \equiv (dx_{1}, dx_{2}, dx_{3})$.

At this point there is still some gauge freedom available to further
specify the connection $A$, namely that provided by 
``symmetric'' gauge transformations. Such a transformation $g$ 
is determined by 
its values on the slice $L$, which must be 
$SO(2)$-invariant.  This yields $i\,g(x,0,0)=g(x,0,0)\,i,$ thus $g(x, 0,
0)= exp\,(i\,h(x)).$ Therefore, symmetric gauge transformations are of
the type
\begin{equation}    \label{symg}
    g(x_{1}, x_{2}, x_{3}) = 
    exp\,(h(r)\,\frac{\sigmadotex}{\vert\overex\vert})\;.
    \end{equation}
where $h(r)$ is an arbitrary function of the radius $r$.
Performing such a gauge transformation does not change the 
form \eqref{symphi} of the Higgs field $\phi.$ 
However, $A$ transforms nontrivially. 
In particular, setting $h(r)\equiv \int \frac{\alpha (r)}{r} \,dr,$ 
exactly cancels the $\alpha$-piece in $g^{-1}dg + g^{-1}A\,g.$
Thus, one can impose $\alpha(r)=0$. The final version of $A$ is then 
    \begin{equation}
	\label{symA}
	 A( x_{1}, x_{2}, x_{3})= 
 \frac{\beta(r)}{r^{3}}(\excrossdex)\cdot (\excrosssigma) 
 +\frac{\gamma(r)}{r^{2}}\oversigma \cdot (\overex \times \overdex)\,.  
 \end{equation}
The only gauge freedom remaining is from the constant of integration in
the indefinite integral $\int  \frac{\alpha (r)}{r} \,dr$ (cf. Section 4).

Note that the $\beta$ and $\gamma$ terms have opposite parities.  The
isometry $\vec x \to -\vec x$ sends $\beta$ to $-\beta$ but sends
$\gamma$ to $+\gamma$.

\section{The Yang-Mills-Higgs functional on spherically 
symmetric configurations}
\setcounter{equation}{0}
In this section we explicity compute the Yang-Mills-Higgs functional,
and the resulting equations of motion, for symmetric pairs $(A,\phi)$.
Our connection is a sum of two terms,
$$A= B + C\,,$$
where
\begin{equation}
    \begin{split}
B&= \frac{\beta(r)}{r^{3}}(\excrossdex)\cdot (\excrosssigma) 
=\frac{\beta(r)}{r^{3}}
[(r^{2} - x^{2}_{1}) i - x_{1}x_{3} k - x_{1} x_{2} j] d\,x_{1} + \\
&+\frac{\beta(r)}{r^{3}}
[(r^{2} - x^{2}_{2}) j - x_{2}x_{3} k - x_{1} x_{2} i] d\,x_{2} 
+\frac{\beta(r)}{r^{3}}
[(r^{2} - x^{2}_{3}) k - x_{1}x_{3} i - x_{2} x_{3} j] d\,x_{3} 
\,,\quad 
\end{split}
\end{equation}
and 
\begin{equation}
    \begin{split}
\qquad C&= \frac{\gamma(r)}{r^{2}} \oversigma \cdot (\overex \times \overdex)=
\frac{\gamma(r)}{r^{2}}( x_{3}j- x_{2} k) d\,x_{1} + \\
& +
\frac{\gamma(r)}{r^{2}}( x_{1}k- x_{3} i) d\,x_{2} + 
\frac{\gamma(r)}{r^{2}}( x_{2}i- x_{1} j) d\,x_{3} 
\end{split}
\end{equation}

\noindent{\bf Preliminary computations}

The curvature of $A$ is given by 
 $$F= dA + A\wedge A= dB + dC + B\wedge B + C\wedge C\,,$$
 since the cross term
 $B\wedge C + C\wedge B$ is identically zero. 
 
Computing $F_{12}$, the various nonzero terms are 
 \begin{equation}
 (d\,B)_{12} =
\frac{\beta^{\prime}} {r^{2}}(- x_{2}i + x_{1} j)
 \end{equation}
 \begin{equation}
 (d\,C)_{12} =
{(\frac{-\gamma}{r^{2}})}^\prime\big(\frac{x_{1}x_{3}}{r}\,i + 
 \frac{x_{2}x_{3}}{r}\,j\big) 
 +\big(\frac{2\gamma}{r^{2}} + {(\frac{\gamma}{r^{2}})}^\prime
 (\frac{x^{2}_{2}}{r} + \frac{x^{2}_{1}}{r})\big)\,k
 \end{equation}
 \begin{equation}
( B\wedge B)_{12}=\frac{2\beta^{2}}{r^{4}}\,\big[x_{1}x_{3} i + x_{2} x_{3} j
 + x^{2}_{3} k\big] 
 \end{equation}
 \begin{equation}
( C\wedge C)_{12}=
\frac{2{\gamma}^{2}}{r^{4}}\,\big[x_{1}x_{3} i + x_{2} x_{3} j
 + x^{2}_{3} k\big]. 
 \end{equation}
By rotational symmetry, the contributions to $F_{13}$ and $F_{23}$ are 
similar. 

 The covariant derivative of the monopole 
 $\phi = \frac {\varphi(r)}{r} \sigmadotex $ is given by
 \begin{equation}
 D_{A}\phi\equiv d\phi + [B,\phi] + [C, \phi]\,,
 \end{equation}
 where 
 \begin{equation}
     \label{dphi}
 d\phi = \frac {\sigmadotex}{r} \, d\varphi  + 
 \varphi \, d\,\big( \frac {\sigmadotex}{r}\big)\;, 
 \end{equation}
 \begin{equation}
 [B,\phi] \equiv B\phi -\phi B= 
 \frac{2\beta\varphi}{r^{2}} ( -x_{3}j+ x_{2} k) \,d\,x_{1} + \hbox{cyclic
permutations},
 \end{equation}
 \begin{equation}
 [C,\phi]\equiv C\phi -\phi C =
 \frac{2\gamma\varphi}{r^{3}}
 [(r^{2} - x^{2}_{1}) i - x_{1}x_{3} k - x_{1} x_{2} j] d\,x_{1} + 
\hbox{cyclic permutations}.
 \end{equation}

At this point we are ready to compute the three terms involved in the 
Yang-Mills-Higgs functional \eqref{ymh}. They are
\begin{itemize}
\item[1)] $|F|^{2}\equiv F\wedge * F$  
\item[2)] $|D\phi|^{2}\equiv D\phi\wedge * D\phi$ 
\item[3)] $(|\phi|^{2}- 1)^{2}\;.$
\end{itemize} 

{\bf Computation of 1):} One shows easily that 
$$dB\wedge * dC =dC\wedge * dB=0\,,$$
$$dB\wedge * (B\wedge B) = (B\wedge B)\wedge * dB=0\,,$$
$$dB\wedge * (C\wedge C) = (C\wedge C)\wedge * dB=0\,,$$
thus 
\begin{equation}
\begin{split}
|F|^{2} =& |dB|^{2}+ |dC|^{2} + dC \wedge * (B\wedge B)
+ dC \wedge * (C\wedge C) 
+(B\wedge B)\wedge *(B\wedge B) \\
&+(B\wedge B)\wedge *(C\wedge C)
+(C\wedge C)\wedge *(C\wedge C) \\
=& 2\frac{{\beta^{\prime}}^{2}}{r^{2}} + 2\frac{{\gamma^{\prime}}^{2}}{r^{2}} 
+4\frac{\gamma^{2}}{r^{4}} + 4\frac{\beta^{2}\gamma }{r^{4}}
+ 4\frac{\gamma^{3}}{r^{4}} + 4\frac{\beta^{4}}{r^{4}} +
8\frac{\beta^{2}\gamma^{2}}{r^{4}} + 4\frac{\gamma^{4}}{r^{4}} \\
= & \frac{2 \beta'{}^2}{r^2} + \frac {2 \gamma'{}^2}{r^2} + 
\frac {4(\beta^2 + \gamma^2 + \gamma)^2}{r^4}.
\end{split}
\end{equation}


{\bf Computation of 2):} 
One shows that 
$$d\phi\wedge * [B, \phi]= [B, \phi]\wedge * d\phi=0\,,$$
$$[B,\phi]\wedge *[C, \phi]= [C, \phi]\wedge * [B,\phi]=0\,,$$
thus 
\begin{eqnarray}
|D\phi|^{2} & = & |d\phi|^{2} + |[B, \phi]|^{2} +|[C, \phi]|^{2} +
2\,d\phi\wedge * [C, \phi] \\
& = & \left ( {\varphi^{\prime}}^{2}+ 2\frac{\varphi^{2}}{r^{2}}\right)
+ 8\frac{\beta^{2}\varphi^{2}}{r^{2}} + 8\frac{\gamma^{2}\varphi^{2}}{r^{2}}
+ 8\frac{\gamma \varphi^{2}}{r^{2}}.
\end{eqnarray}

{\bf Computation of 3):} One easily obtains 
$$(|\phi|^{2} - 1)^{2}= \big(\frac{\varphi^{2}}{r^{2}}-1\big)^{2}.$$

Collecting terms, 
the Yang-Mills-Higgs functional calculated on spherically symmetric 
configurations is given by  
\begin{equation}\label{symh}
\begin{split}
\SYMH (\gamma, \varphi) 
&= 4 \pi \int_0^1 [\ 2\epsilon \left({\beta^{\prime}}^{2} + {\gamma^{\prime}}^{2}
+
\frac {2}{r^2} (\beta^{2} + \gamma^{2} +\gamma)^{2}\right) +
r^2 {\varphi^{\prime}}^{2} + \\
&\quad 
2 \varphi^2 [1 + 4(\beta^{2} + \gamma^{2} + \gamma)] +\lambda 
r^{2}(\varphi^{2}-1)^{2}\ ] dr \ .
\end{split}
\end{equation}

\section{Further gauge transformations and the Euler-Lagrange 
equations} 
\setcounter{equation}{0}
We want to search for absolute minima of the functional \eqref{symh}
among all finite-action spherically symmetric configurations $(\beta,
\gamma, \varphi).$  The following theorem restricts the possibilities:

\begin{theorem}
If the functional \eqref{symh} has a minimum, then this minimum is
achieved by functions $(\beta, \gamma, \varphi)$ with $\beta$ identically
zero and $\gamma(0)=0$. 
\end{theorem}

\begin{proof}
First we find a gauge transformations that yields
$\beta(0)=\gamma(0)=0$.  In that gauge, we then show that minimization
requires $\frac{\beta}{\gamma + \frac{1}{2}}$ to be constant, hence for $\beta$
to be identically zero. 

If $\SYMH (\beta, \gamma, \varphi) <\infty $, then $\beta$ and
$\gamma$ much approach well-defined limits as $r \to 0$, and 
$\beta^{2}(0) + \gamma^{2}(0) + \gamma(0) = 0$.  If $\beta(0)$ and
$\gamma(0)$ are not already zero, we let $\theta$ be the argument of
the complex number $\beta(0) + i \gamma(0)$, and define
\begin {equation} 
g= exp (\theta 
\frac {\sigmadotex}{\vert \overex\vert })
\equiv  exp (\theta \, 
\frac{x_{1} i + x_{2} j + x_{3} k} {\vert \overex\vert})
\equiv  \cos(\theta) + \sin(\theta)\, 
\frac {\sigmadotex }{\vert \overex\vert}\  .
\end{equation}
Then  
\begin {equation}
d\,g (x_{1}, x_{2}, x_{3})= \sin(\theta) \left [
\frac {\sigmadotdex}
{\vert \overex\vert }
+ d(\frac{1}{\vert \overex\vert} )
\sigmadotex \right ] \;, 
\end{equation}
and evaluating on our slice gives 
\begin{equation}
d\,g(r,0,0) = \frac{\sin(\theta)}{r} (j\,d\,x_{2} + k\, d\,x_{3}).
\end{equation}
Our transformed connection on the slice is then 
\begin{equation}
    \begin{split}
[g^{-1}d\,g+ g^{-1}A\,g] (r,0,0) = &
 \left[ \frac {\cos (2\theta)\,\beta(r)}{r} + 
\frac {\sin (2\theta)\,\gamma(r)}{r} 
+\frac{\cos(\theta) \sin(\theta)}{r}\right]
(j\,d\,x_{2} + k\, d\,x_{3})\\
&+
\left[-\frac{\sin (2\theta) \,\beta(r)}{r} 
+ \frac {\cos (2\theta)\,\gamma(r)}{r} 
-\frac{\sin^{2}(\theta)}{r}\right]
(k\,d\,x_{2} -j\, d\,x_{3}) \\
\equiv & \frac{\hat\beta(r)}{r}
(j\,d\,x_{2} + k\, d\,x_{3})
+\frac{\hat\gamma (r)}{r} 
(k\,d\,x_{2} - j\, d\,x_{3}).
\end{split} 
\end{equation}
Plugging in the values of $\sin(\theta)$, $\cos(\theta)$, etc., gives
\begin{eqnarray}
\hat \beta(r) & = & \frac{(\beta^2(0)-\gamma^2(0)) \beta(r) + 
2\beta(0)\gamma(0) \gamma(r) + \beta(0) \gamma(0)}{\beta^2(0)+\gamma^2(0)}, \\
\hat \gamma(r) & = & \frac{(-2\beta(0)\gamma(0)) \beta(r) + 
(\beta^2(0)-\gamma^2(0)) \gamma(r) - \gamma^2(0)}{\beta^2(0)+\gamma^2(0)}.
\end{eqnarray}
As $r \to 0$, both these terms go to zero, since $\beta^2(0)+\gamma^2(0)
+ \gamma(0)=0$.

Having set $\beta(0)=\gamma(0)=0$, we now show that $\beta$ must be
identically zero.  We choose polar coordinates in the 
$(\beta, \gamma+\frac{1}{2})-$plane:
\begin{equation}\label{betagamma}
\left\{ \begin{split}
&\beta= \nu \,\cos\,t\\
&\gamma+\frac{1}{2}=\nu \sin\, t
\end{split} \right.
\end{equation} 
In these coordinates, the functional
\eqref {symh} becomes

\begin{equation}\label{ymhsympol}
\begin{split}
\SYMH (\nu, t) 
&= 4 \pi \int_0^1 [\ 2\epsilon \left({\nu^{\prime}}^{2} + 
{t^{\prime}}^{2}\nu^{2}
+
\frac {2}{r^2} (\nu^{2} -\frac{1}{4})^{2}\right) +
r^2 {\varphi^{\prime}}^{2}+ \\
&\quad 
8 \varphi^2\nu^{2} +\lambda 
r^{2}(\varphi^{2}-1)^{2}\ ] dr \ ,
\end{split}
\end{equation}
with $\nu ^{2}(0)=\frac{1}{4}.$
The only dependence on $t$ is in the ${t'}^2 \nu^2$ term, which is minimized
by setting $t=$ constant. But then $\cot(t) = 
\frac{\beta}{\gamma + \frac{1}{2}}$ must be constant, and equal to
$\frac{\beta(0)}{\gamma(0) + \frac{1}{2}}=0$, so $\beta$ is identically zero.
\end{proof}

We may therefore restrict our attention to 
the functional 
\begin{equation}\label{symhgamma}
{\mathcal S} (\gamma, \varphi) 
= 4 \pi \int_0^1 [\ 2\epsilon \left({\gamma^{\prime}}^{2} +
\frac {2}{r^2} (\gamma^{2} +\gamma)^{2}\right) +
r^2 {\varphi^{\prime}}^{2} +
2 \varphi^2  (1+2\gamma)^{2} +\lambda 
r^{2}(\varphi^{2}-1)^{2}\ ] dr \ .
\end{equation}
The Euler-Lagrange equations for this functional together with the 
appropriate boundary conditions are then
\begin{equation}
    \label{ode}   
\left\{ \begin{split}
&\gamma^{\prime\prime} -\frac{2}{\epsilon} \varphi^{2}(1+2\gamma) 
-\frac{2}{r^{2}} (\gamma^{2}+ \gamma) (1+2\gamma) =0  \qquad\text { on } 
(0,1)\\
&\varphi^{\prime\prime}  + \frac{2\varphi^{\prime}}{r} -
\frac{2\varphi}{r^{2}} \,(1+2\gamma)^{2} - 
2\lambda\varphi(\varphi^{2}-1) = 0 \qquad\text { on } (0,1)\\ 
&\gamma (1)=-\frac{1}{2}  \\  
&\gamma (0) = 0 \\
&\varphi (1) = + 1 \quad(\text { or } \varphi (1) = -1) 
\end{split}\right.
    \end{equation} 
The boundary conditions above come directly from the variational principle.
In fact, to cancel the boundary terms one needs to either restrict 
the space of connections to those with prescribed boundary data, or to impose 
Neumann-type boundary conditions.
The boundary condition $\varphi(1) = \pm 1$ comes from $|\phi|=1$ on
$\partial B^3$ and $\gamma(1) = -\frac{1}{2}$ comes from
\begin{equation}
(D \phi)_\tau = (d\phi)_\tau + [A_\tau, \phi] = (1 + 2\gamma) \phi
\; d \! \left (\frac{\vec \sigma \cdot \vec x}{r} \right ) = 0
\end{equation}
on $\partial B^3$, where the subscript $\tau$ denotes tangential components
(cf. \cite {mar1, mar2}).

{\bf Observation:} After some computation, the system \eqref{ode}
could also be obtained by imposing spherical symmetry in \eqref{pde},
thus critical symmetric points for the action \eqref{symhgamma} are
symmetric critical points (not necessarily minima) for \eqref{ymh}.
This is also known apriori from the ``principle of symmetric
criticality'' \cite{pal}.

\section{Existence of spherically symmetric monopoles}
\setcounter{equation}{0}

Our basic existence result is
\begin{theorem}
\label{existence}
For all values of $\lambda\geq 0,$ $\epsilon > 0,$
there exists a symmetric solution of 
\begin{equation} \label{system}
   \left\{ \begin{split}
&\epsilon\,*D_A* F=[D_A\phi,\,\phi] \qquad\text { on } B^{3}\\
&*D_A*D_A\phi=\frac{\lambda}{2}(|\phi|^2-1)\phi  \qquad\text { on } B^{3}\\
&(D\phi)_{\tau} = 0   \qquad\text { on }  \partial B^{3} \\
&\vert \varphi\vert = 1  \qquad\text { on }  \partial B^{3} 
\end{split}\right.
\end{equation} 
\end{theorem}
      {\bf Observation:}
These equations do not reduce to the Bogomolnyi equations, even if
$\lambda= 0.$  The Bogolmolnyi argument involves an integration by parts;
on a finite domain, this results in a boundary contribution \cite{ja-ta}.

\begin{proof}  Because of the derivative terms in the action, the natural
space for $\gamma$ (denoted ${\mathcal H}_\gamma$) is $H^1(0,1)$, while the 
natural space ${\mathcal H}_\varphi$ for $\varphi$ is
the weighted Sobolev space $H^1((0,1),r^2 dr)$.
By the Sobolev embedding theorem, functions in ${\mathcal H}_\gamma$ 
are continuous on $[0,1]$.  Functions in ${\mathcal H}_\varphi$ are
continuous on $(0,1]$, but may not have a limit at $r=0$.  
We may therefore apply boundary conditions to
$\gamma$ at $r=0$ and at $r=1$, and to $\varphi$ at $r=1$.
 
Let 
\begin{equation}
{\mathcal F}=\{(\gamma, \varphi)\in {\mathcal H}_\gamma \times {
\mathcal H}_\varphi \;:\; 
\gamma(1)=-\frac{1}{2},\,\gamma(0)=0,
    \, \varphi (1) = 1. \}
\end{equation}
%
%
The action functional \eqref{symhgamma} is well-defined on ${\mathcal
F}$, and is finite whenever $\varphi$ is bounded. In particular, $\mu
\equiv Inf_{\mathcal F}\,{\mathcal S}$ is finite.  We follow the
direct method in the calculus of variations.  That is, take a
minimizing sequence for ${\mathcal S}$, show that it converges weakly
in ${\mathcal F}$, and then show that the weak limit minimizes the
action and so solves the Euler-Lagrange equations.

Let $(\gamma_n,\varphi_n)$ be a minimizing sequence for ${\mathcal
S}$.  Since $\lambda \ge 0$, the action is not increased if we 
make the replacement
\begin{equation} 
\varphi(r) \to 
\begin{cases}
 -1, & \hbox{if }\varphi(r) < -1; \cr
\varphi(r), & \hbox{if }-1 \le \varphi(r) \le 1; \cr
1, & \hbox{if }\varphi(r) > 1.
\end{cases}
\end{equation}
As a result, we can assume that each $\varphi_n(r)$ is bounded in
magnitude by 1.  Under these circumstances, the sequence $(\gamma_n,
\varphi_n)$ is bounded in ${\mathcal F} \subset {\mathcal H}_\gamma
\times {\mathcal H}_\varphi$.  However, balls in ${\mathcal H}_\gamma$
are weakly compact, as are balls in ${\mathcal H}_\varphi$, so the
pair $(\gamma_n, \varphi_n)$ converges weakly in ${\mathcal H}_\gamma
\times {\mathcal H}_\varphi$ to a limit $(\gamma_\infty,
\varphi_\infty)$.

By Sobolev, $\gamma_n(r)$ and $\varphi_n(r)$ converge pointwise to
$\gamma_\infty(r)$ and $\varphi_\infty(r)$, so the limiting values
$\gamma(0)$, $\gamma(1)$, and $\varphi(1)$ are preserved, and
$(\gamma_\infty, \varphi_\infty) \in {\mathcal F}$. Moreover, terms in
${\mathcal S}(\gamma_n,\varphi_n)$ that don't involve derivatives
converge to the corresponding terms in ${\mathcal S}
(\gamma_\infty,\varphi_\infty)$.  The derivative terms are
quadratic, hence weakly semicontinuous.  As a result, ${\mathcal
S}(\gamma_\infty, \varphi_\infty)$ is bounded above by $\mu$, and
therefore must equal $\mu$.

Showing that $\gamma_\infty$ and $\varphi_\infty$ satisfy the
Euler-Lagrange equations \eqref{ode} is then a standard exercise in
the calculus of variations.  Smoothness of $(\gamma_\infty,
\varphi_\infty)$ away from $r=0$ follows by elliptic regularity of the
equations \eqref{ode}.  Smoothness at $r=0$ follows from regular
singular-point analysis, combined with the fact that both functions
are bounded (see \S 6 for details).  
This in turn implies that the
connection and Higgs field $(A,\phi)$ constructed from $(\gamma_\infty,
\varphi_\infty)$ comprise a smooth, symmetric classical solution to
the PDE system \eqref{system}.  (Alternatively, one can establish
regularity of $(A,\phi)$ from the ellipticity of the PDE system 
\eqref{system}, since we are
working in a gauge with $d^* A=0$.)
 
\end{proof}

\section{Regular Singular Point Analysis}

In \S 5 we demonstrated the existence of symmetric monopoles for arbitrary
$\lambda \ge 0$ and $\epsilon > 0$.  In this section we explore their 
form near the regular singular point of the equations \eqref{ode}, 
namely $r=0$. 

\begin{theorem}
Let $(\gamma, \varphi)$ be a bounded finite-action 
solution to the ODE system \eqref{ode}
for some fixed $\epsilon >0$ and $\lambda \ge 0$.  Then there exist
constants $a_1$ and $b_2$ such that 
\begin{equation}
\begin{split}
\varphi(r) = & a_1 r + O(r^3); \\
\gamma(r) = & b_2 r^2 + O(r^4); \\
\gamma'(r) = & 2 b_2 r + O(r^3);
\end{split}
\end{equation}
near $r=0$. In particular, $\varphi(0)=\gamma'(0)=0$.
\end{theorem}

\begin{proof}
We begin with the equation for $\varphi$, which we write as
\begin{equation}
\label{phieq}
\varphi'' + \frac{2\varphi'}{r} - \frac{2 \varphi}{r^2}
= \frac{8 \varphi}{r^2} (\gamma + \gamma^2) + 2 \lambda \varphi (\varphi^2-1).
\end{equation}
Since $\varphi$ is bounded and $\gamma(0)=0$, the terms on the right hand side
are less singular than those on the left hand side, and to leading order
$\varphi$ resembles the solution to the homogeneous equation
\begin{equation}
\varphi'' + \frac{2\varphi'}{r} - \frac{2 \varphi}{r^2} = 0.
\end{equation}
The general solution to this equation is $\varphi = a_1 r + a_{-2}r^{-2}$.
However, $\varphi$ is bounded by assumption, so $a_{-2}=0$. 
Thus the solution to \eqref{phieq}
is, to leading order,  $a_1 r$.

Next we turn to the equation for $\gamma$, namely
\begin{equation}
\label{gammaeq}
\gamma'' - \frac{2 \gamma}{r^2} = \frac{2}{\epsilon} \varphi^2 (1+2\gamma)
+ \frac{2\gamma}{r^2}(2 \gamma^2 + 3\gamma).  
\end{equation}
Again, since $\varphi$ is bounded and $\gamma(0)=0$, this may be viewed as 
a perturbation of the homogeneous linear equation
\begin{equation}
\gamma'' - \frac{2 \gamma}{r^2} = 0,
\end{equation}
whose solution is $\gamma = b_{-1} r^{-1} + b_2 r^2$.  Since $\gamma$ is 
bounded, $b_{-1}=0$. Thus our solution to \eqref{gammaeq} is, to leading
order, $b_2 r^2$. 

With these basic results, we can estimate the right hand sides of 
\eqref{phieq} and \eqref{gammaeq}.  The right hand side of \eqref{phieq} is
$O(r)$, which gives an $O(r)$ correction to $\varphi''$, hence an $O(r^3)$
correction to $\varphi$.  The right hand side of \eqref{gammaeq} is
$O(r^2)$, thus giving an $O(r^3)$ correction to $\gamma'$ and an
$O(r^4)$ correction to $\gamma$. 
\end{proof} 

One can do an expansion for $\varphi$ and $\gamma$ in powers of $r$.
Indeed, only odd powers contribute to $\varphi$ and only even powers
contribute to $\gamma$.  This is seen by induction.  By Theorem 6.1, 
$\varphi$ is odd and $\gamma$ is even through order $r^2$.  However, if
$\varphi$ is odd and $\gamma$ is even through order $r^k$, then the
right hand sides of \eqref{phieq} and \eqref{gammaeq} are odd and even,
respectively, to order $r^k$, and so $\varphi$ and $\gamma$ are odd 
and even, respectively, to order $r^{k+2}$. Thus $\varphi$ and $\gamma$
are odd and even to all orders in $r$. 

We can therefore write an asymptotic expansion
\begin{equation} \label{powers}
\begin{split}
\varphi(r) & \sim \sum_{n\hbox{\small{} odd}} a_n r^n, \\
\gamma(r) & \sim \sum_{n\hbox{\small{} even}} b_n r^n.
\end{split}
\end{equation}
Plugging this expansion into equations \eqref{phieq} and \eqref{gammaeq}
and equating coefficients of $r^{n-2}$ gives recursion relations of the form
\begin{equation}\label{recursion}
\begin{split}
(2k)(2k+3) a_{2k+1} = & \hbox{ algebraic expression involving }
a_1, b_2, \ldots,b_{2k}, \\
(2k+1)(2k-2) b_{2k} = & \hbox{ algebraic expression involving }
a_1, b_2, \ldots,b_{2k-2}. 
\end{split}
\end{equation}
These relations do not constrain $a_1$ or $b_2$.
However, once we have $a_1$ and $b_2$, the remaining coefficients
are determined.  The first few are
\begin{equation}
\begin{split}
a_3 = & (4 a_1 b_2 -  \lambda a_1)/5; \\
b_4 = & \left ( 3 b_2^2 + \epsilon^{-1} a_1^2 \right )/5; \\
a_5 = & (4 a_1 b_4 + 4 a_3 b_2 + 4 a_1 b_2^2 + \lambda(a_1^3 - a_3))/14; \\
b_6 = & \left (  b_2^3 + 3 b_2 b_4 + \epsilon^{-1}
(a_1 a_3 + a_1^2 b_2)
\right ) / 7; \\
a_7 = & [4(  a_1 b_6 +  a_3 b_4 +  a_5 b_2  + a_3 b_2^2 + 2 a_1 b_2 b_4)
+ \lambda(3 a_1^2 a_3 - a_5)]/27; \\
b_8 = & \left [ 3(2 b_2^2b_4 + b_4^2 +  2 b_2 b_6) + \epsilon^{-1}
(a_3^2 + 2 a_1 a_5
+ 2 a_1^2 b_4 + 4 a_1 a_3 b_2 ) \right] / 27; \\
a_9 = & \big [4(a_1 b_8 + a_3 b_6 + a_5 b_4  + a_7 b_2  + a_5 b_2^2  
 + 2 a_3 b_2 b_4 + a_1 b_4^2 + 2 a_1 b_2 b_6) \\
& + \lambda(3 a_1^2 a_5 + 3 a_1 a_3^2 - a_7)\big ]/44; \\
b_{10} = & \big [ 3(b_2^2 b_6 + b_2 b_4^2 + b_2 b_8 + b_4 b_6) \\
& + \epsilon^{-1} (a_1 a_7 + a_3 a_5 + a_1^2 b_6 + a_3^2 b_2 
+ 2 a_1a_3b_4 + 2 a_1a_5b_2)
\big ] / 22. 
\end{split}
\end{equation}

\section{Symmetries and Stability}


The action functional and the resulting Euler-Lagrange equations 
are invariant under two natural symmetries.
\begin{eqnarray}
\label{phisym}
\varphi(r) \to - \varphi(r) & \qquad & \gamma(r) \to + \gamma(r); \\
\label{gammasym}
\varphi(r) \to \varphi(r) & \qquad & \gamma(r) \to -1 - \gamma(r).
\end{eqnarray}
The first symmetry \eqref{phisym}
comes from the isometry $\vec x \to - \vec x$ of $B^3$,
which of course respects rotational symmetry.  
Since $\vec \sigma \cdot \vec x$ is odd and $\vec \sigma \cdot (\vec x \times
\vec {dx})$ is even, pulling the pair $(A,\phi)$ back by this isometry flips
the sign of $\varphi$ while preserving $\gamma$.  Using this symmetry, we
can fix the sign of $\varphi(1)$, which we henceforth take to be positive.

The second symmetry \eqref{gammasym}
is a gauge transformation by $(\vec\sigma \cdot \vec x)/r$.
This is of the form of (4.2), with $\theta = \pi/2$.  From (4.4) it is 
clear that this transformation sends $\gamma$ to $-1-\gamma$ without generating
a $\beta$ term or changing $\varphi$.  Applying this to a connection with
$\gamma(0)=0$ yields a new connection with $\gamma(0)=-1$.  This connection
has finite action but is singular at the origin, reflecting the singular
gauge transformation that generated it.  

We now consider stability properties of the ODE system \eqref{ode}.  
These ODEs have several fixed points, namely
\begin{equation}
(\gamma,\varphi) =  \hbox{ $( -\frac{1}{2}, 1 )$, $(-\frac{1}{2}, -1)$, 
$(-\frac{1}{2}, 0 )$, $(0,0)$  or  $(-1,0)$}.
\end{equation}
The point $\left (-\frac{1}{2}, -1 \right )$ is related to 
$\left ( -\frac{1}{2}, 1 \right)$ by the symmetry \eqref{phisym}, while
$(-1,0)$ is related to $(0,0)$ by \eqref{gammasym}, so we do not need to
study these.  What remains is  $\left ( -\frac{1}{2}, 1 \right)$,  
$\left (-\frac{1}{2}, 0 \right)$, and (0,0).

For the $\gamma = -1/2$ fixed points, we define $\delta = \gamma + 1/2$,
and the equation for $\gamma$ becomes
\begin{equation}
\delta'' = 4 \delta \left ( \frac{\varphi^2}{\epsilon}
- \frac{1}{4 r^2} + \frac{\delta^2}{r^2} \right ).
\end{equation}
The fixed point $\gamma = -1/2$ is stable for $\varphi=1$ when $r^2 <
\epsilon/4$, but is unstable if $r^2 > \epsilon/4$.  This defines a
natural length scale to the problem, namely $\sqrt{\epsilon}/2$.  We
should expect our solutions to behave qualitatively differently for
$r$ less than or greater than this length scale.  Of course, if
$\epsilon>4$, then all radii $r$ are less than this length scale.  In
the case of $\varphi=0$, the value $\gamma=1/2$ is always stable,
regardless of $\epsilon$ or $r$.

Near $\gamma=-1/2$, the equation for $\varphi$ becomes
\begin{equation}
\varphi'' + 2 \frac{\varphi'}{r} = 2 \varphi \left (
\frac{4 \delta^2}{r^2} + \lambda(\varphi^2-1) \right ).
\end{equation}
The behavior of this fixed point depends on the value of $\varphi$. 
Near $\varphi=0$ we have
$$ (r\varphi)'' = -2\lambda r\varphi + O(\delta^2) + O(\varphi^2),$$
which is stable for all positive values of $r$. 
Near $\varphi=1$, however, we write $\varphi = 1 + \zeta$ and have
$$ (r\zeta)'' = 4 \lambda (r \zeta) + O(\zeta^2)
+ O(\delta^2).
$$
This is unstable as long as $\lambda >0$, and has natural length
scale $1/\sqrt{4\lambda}$. 

To summarize, the fixed point $(-1/2,0)$ is stable, while the fixed point
$(-1/2,1)$ is unstable.  If $r^2 < \epsilon/4$, then there is one unstable
mode, corresponding to growth of $\varphi-1$.  If $r^2 > \epsilon/4$, then
there are two unstable modes, one for $\varphi$ and one for $\gamma$. 

Finally, there is the fixed point (0,0).  Near (0,0) we have
\begin{equation}
\begin{split}
 \gamma'' = & 2 \gamma/r^2 + \hbox{ higher order}, \\
\varphi'' + 2 \varphi'/r = & 2 
\varphi (r^{-2} - \lambda) +\hbox{ higher order.} 
\end{split}
\end{equation}
This fixed point is always unstable, with $\gamma$ growing rapidly. $\varphi$
will grow exponentially if $r < 1/\sqrt{\lambda}$, and will oscillate if 
$r > 1/\sqrt{\lambda}$.

\vbox{
\centerline{\epsfysize=4truein\epsfbox{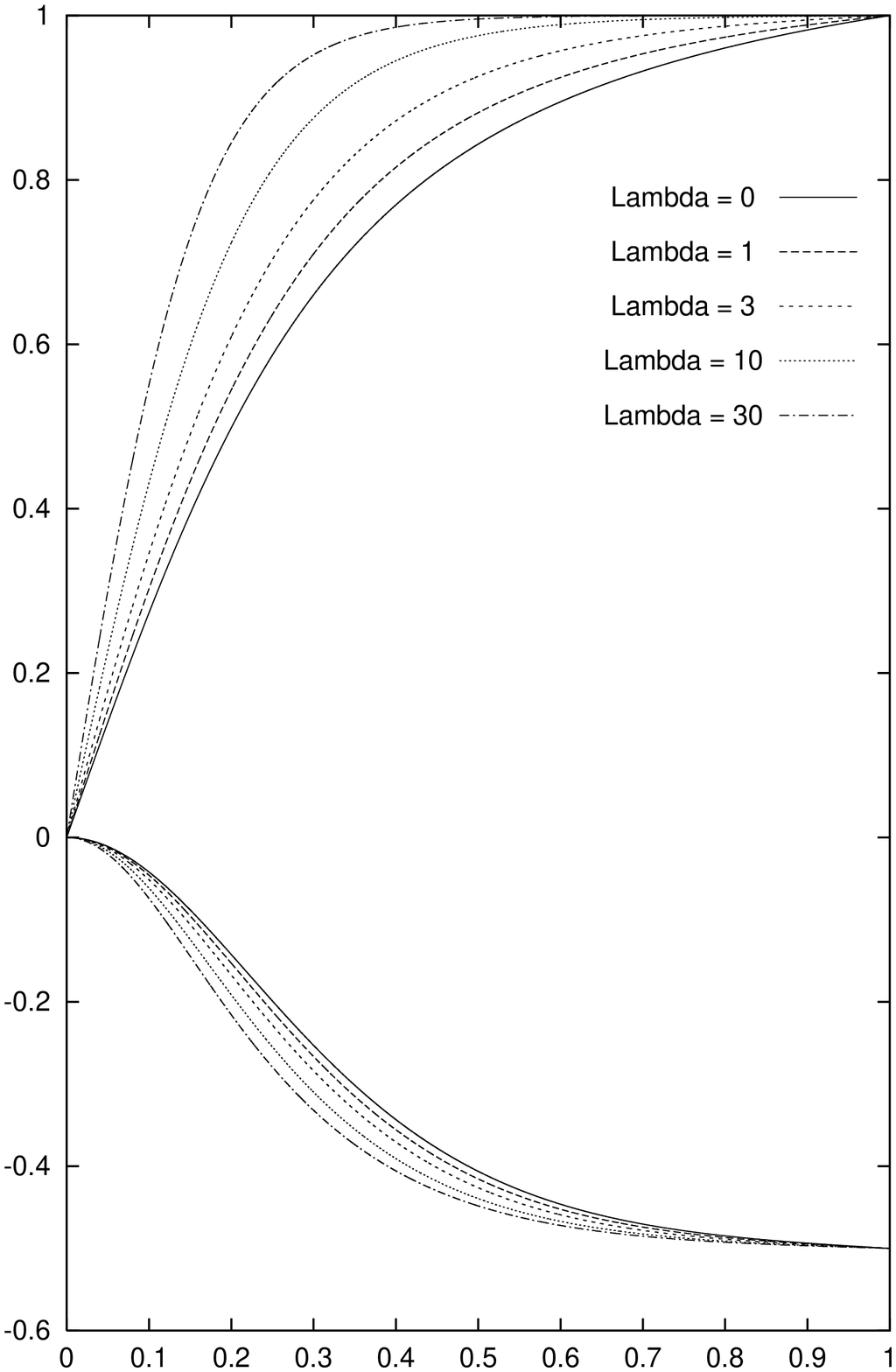}, 
\epsfysize=4truein\epsfbox{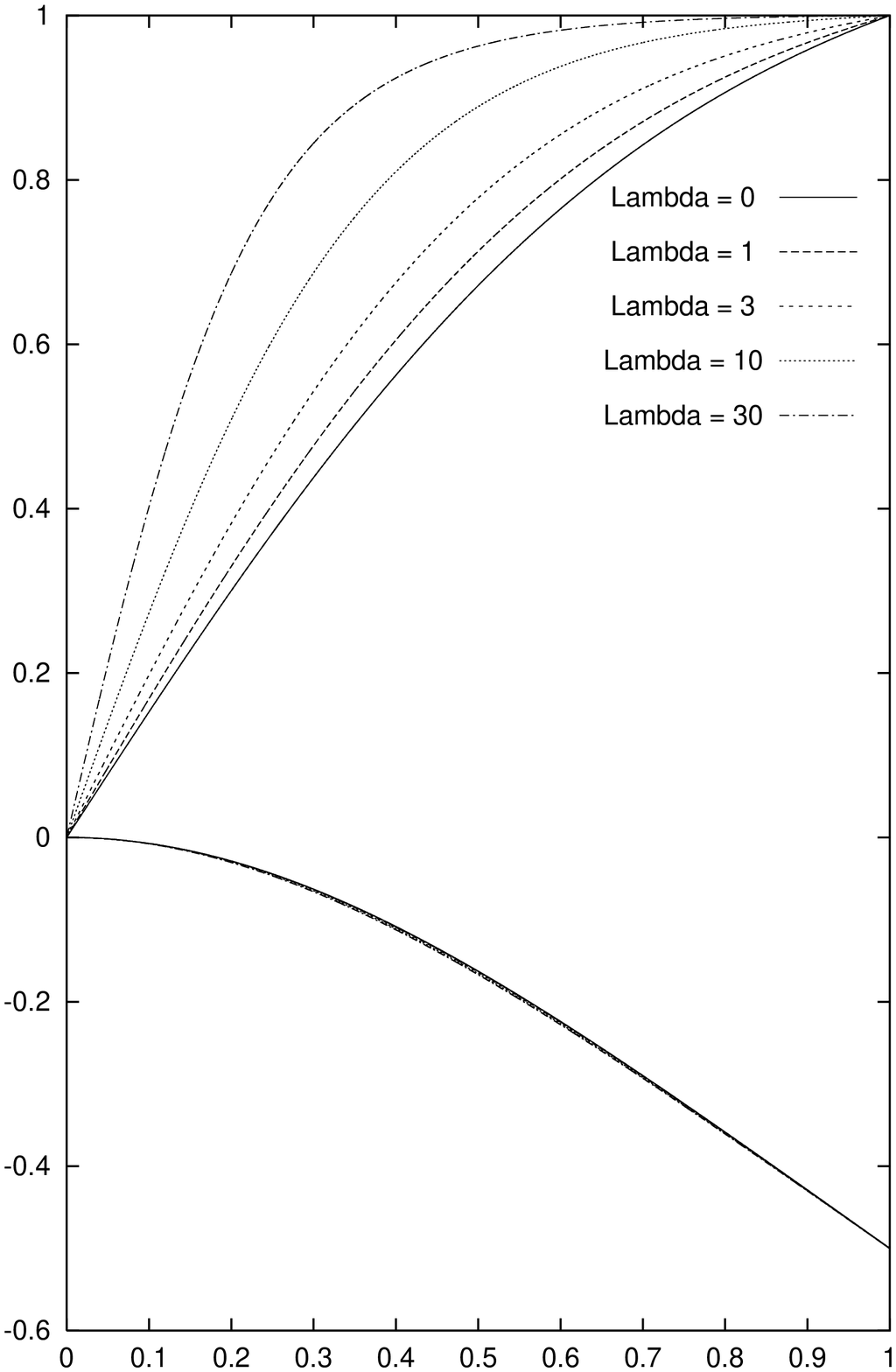}}
\smallskip
\centerline{\bf $\epsilon=0.1$ \hskip 2 in $\epsilon = 10$}
\smallskip
\centerline{\bf Figure 1. Trajectories with fixed $\epsilon$.}
\bigskip}

\section{Numerical Results and Qualitative Analysis}

For any fixed $\epsilon$ and $\lambda$, and 
given $a_1$ and $b_2$, one can in principle integrate the
differential equations out to $r=1$.  In practice, numerical errors due to
the discretization of the interval $[0,1]$ can be very bad near the origin,
due to the singular nature of the ODE system there.  A better method is
to use the power series \eqref{powers} in a neighborhood of the origin
and to numerically integrate from there.  In a discretization of $10,000$ 
points, we use the power series out to $r=0.01$, or 100 lattice spacings
from the origin.  

In this way we get a pair $(\gamma(1), \varphi(1))$ for each $(a_1, b_2)$.
Using Newton's method, we then find values of $(a_1,b_2)$ such that 
$(\gamma(1), \varphi(1))=(-1/2, 1)$.  Table 1 lists the correct values
of $a_1$ and $b_2$ for several values of $\epsilon$ and $\lambda$.

The resulting functions $\varphi(r)$ and $\gamma(r)$ are sketched in
Figures 1 and 2. Figure 1 shows the functions for different values of
$\lambda$ and $\epsilon$ fixed at 0.1 or at 10.  Figure 2 is similar, only 
with $\lambda$ fixed and $\epsilon$ variable. 
In each case the positive function is $\varphi$
and the negative function is $\gamma$.

\vbox{
\centerline{\epsfysize=4truein\epsfbox{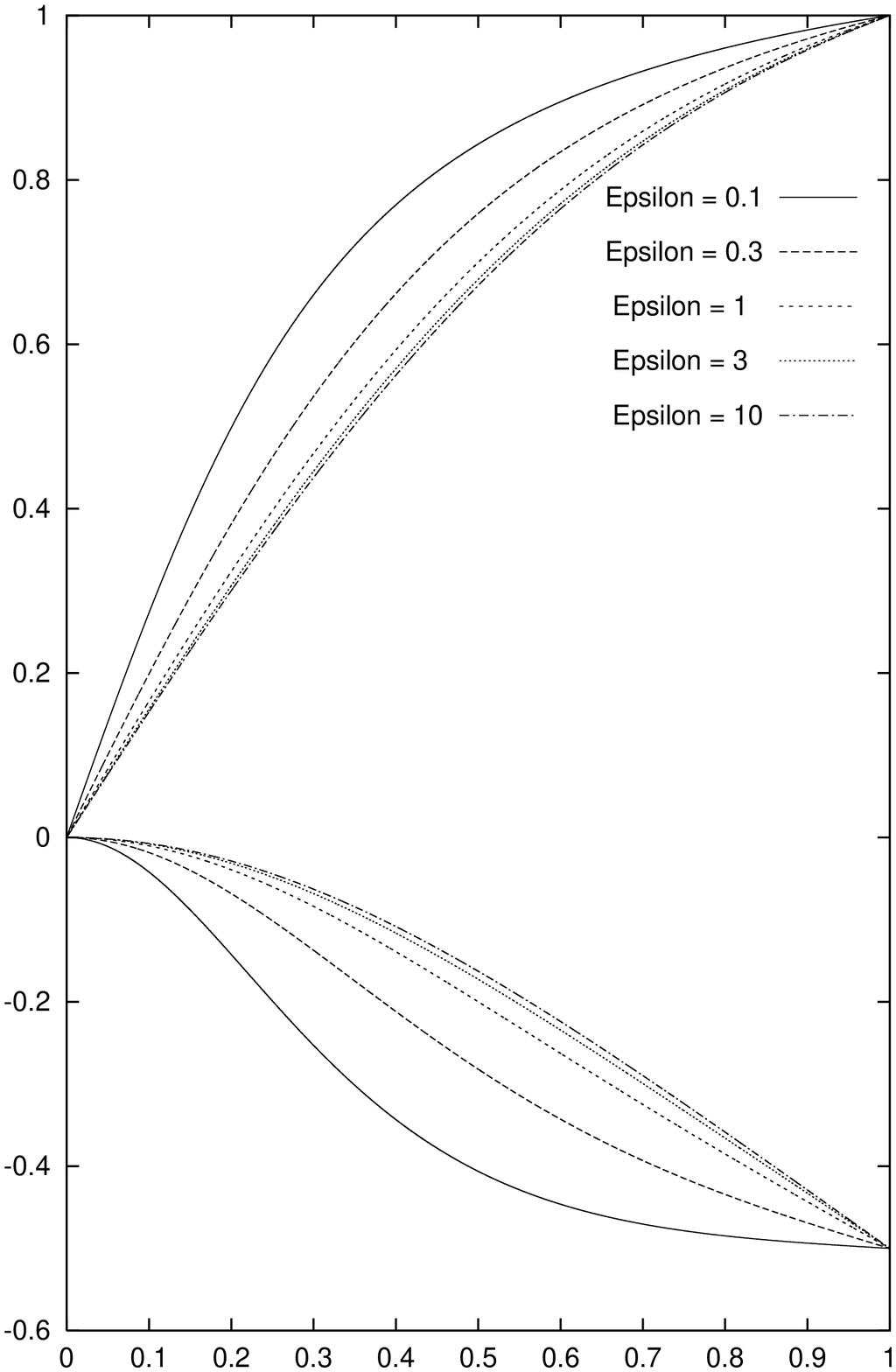}, \epsfysize=4truein\epsfbox{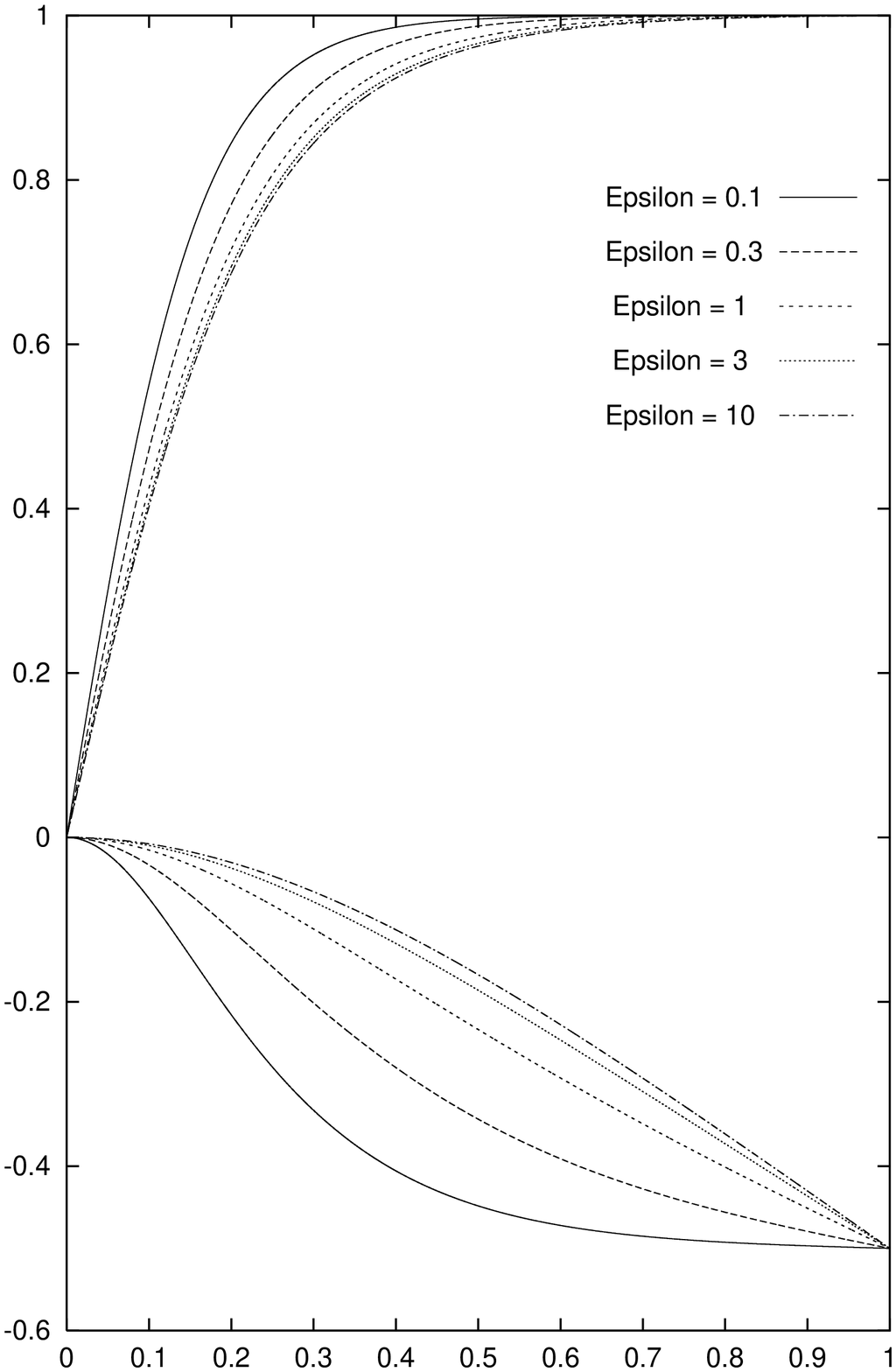}}
\smallskip
\centerline{\bf $\lambda=0$ \hskip 2 in $\lambda = 30$.}
\smallskip
\centerline{\bf Figure 2: Trajectories with $\lambda$ fixed.}
\bigskip}

\begin{table}\vbox{
\begin{tabular}{|r|r|r|r|} \hline
\hfil $\epsilon$ \hfil & \hfil $\lambda$ \hfil &  \hfil $a_1$ \hfil &
\hfil $b_2$ \hfil \\ \hline
0.1 &  0 &  2.82909077 &  -4.47460232 \\ \hline  
0.1 & 1 &   3.14773551 &  -4.92072556 \\ \hline 
0.1 &  3 &  3.62692766  & -5.57110938 \\ \hline 
0.1 &  10 &  4.62892407  & -6.81947999 \\ \hline 
0.1 &  30 & 6.19274693 &  -8.46474894 \\ \hline
0.3 &  0 &  2.01904955 &  -1.88549902 \\   \hline
0.3 & 1 &  2.26118176 &  -2.04994984  \\   \hline
0.3 & 3 &  2.66517994 &  -2.31673622  \\   \hline
0.3 &  10 &   3.59550462 &  -2.86817001  \\   \hline
0.3 & 30 &  5.12510342 &  -3.58374045 \\   \hline
1 & 0 &  1.67098122 &  -1.02894746  \\   \hline
1 & 1 & 1.85973704  &  -1.07504639  \\   \hline
1 & 3 & 2.19572981  &  -1.15577783  \\   \hline
1 & 10 & 3.04898441 &  -1.34041824  \\   \hline
1 & 30 & 4.55384341 &  -1.58910470 \\   \hline
3 & 0 &  1.57081044 &  -0.80615986 \\   \hline
3 & 1 &  1.74184236 &  -0.82078859 \\   \hline
3 & 3 &  2.05156143 &  -0.84695210 \\   \hline
3 & 10 &  2.86750186 &  -0.90924487  \\   \hline
3 & 30 &  4.35776101 &  -0.99551235 \\   \hline
10 & 0 &  1.53622287 &  -0.73146686 \\   \hline
10 & 1 &   1.70099654 &  -0.73576432 \\   \hline
10 & 3 &    2.00102288 &   -0.74350147 \\   \hline
10 & 10 &  2.80198139  & -0.76219107  \\   \hline
10 &  30 &  4.28571713 &  -0.78838676 \\   \hline
\end{tabular}}
\bigskip
\caption{Taylor coefficients 
$(a_1,b_2)$ for various values of $(\epsilon, \lambda)$.}
\end{table}


From these figures several qualitative features are clear.  Although
$\varphi$ depends significantly on both $\epsilon$ and $\lambda$,
$\gamma$ is practically independent of $\lambda$, especially when
$\epsilon$ is large.  The length scale on which $\gamma$ changes from
0 to $-1/2$ is the smaller of $\sqrt{\epsilon}$ and 1.  The length
scale on which $\varphi$ changes from 0 to 1 is the smallest of
$\sqrt{\epsilon}$, $1/\sqrt{\lambda}$, and 1. Thus changing $\lambda$
has the greatest effect when $\lambda$ is greater than 1, while
changing $\epsilon$ has the greatest effect when $\epsilon<1$.


The source code for these numerical results can be obtained from the
authors.

\section*{Acknowledgements}		
The authors thank Prof. Karen K. Uhlenbeck for useful conversations on 
this subject.

\end{document}